\documentclass{amsart}
\usepackage{amssymb,latexsym}

\newtheorem{theorem}{Theorem}
\newtheorem{proposition}{Proposition}
\newtheorem{corollary}{Corollary}
\newtheorem{definition}{Definition}

\newtheorem{lemma}{Lemma}

\def\proof{\smallskip\noindent {\it Proof: \ }}

\def\endproof{\hfill$\square$\medskip}
\def\Z{\mathbb{Z}}

\newcommand{\lk}{\mbox{\upshape lk}\,}
\newcommand{\skel}{\mbox{\upshape Skel}\,}

\title{A short simplicial $h$-vector and the Upper Bound Theorem}

\author{Patricia Hersh 
\and   Isabella Novik}

\address{Department of Mathematics, Box 354350\\
         University of Washington\\
         Seattle, WA 98195}

\email{hersh@math.washington.edu}
\email{novik@math.washington.edu}
\begin{document}

\begin{abstract}
We verify the Upper Bound Conjecture (UBC) for a class of
odd-dimensional simplicial complexes that in particular includes
all Eulerian simplicial complexes with isolated singularities. 
The proof relies on a new invariant of simplicial complexes ---
a short simplicial $h$-vector.
\end{abstract}

\maketitle
\section{Introduction}
The goal of this note is to prove an extension of
the Upper Bound Theorem for (simplicial) polytopes. 
The main tool in the proof  is a certain new invariant
of simplicial complexes, which is a simplicial analog of a short
cubical $h$-vector introduced by Adin \cite{Adin}.

We start by recalling several definitions.
A (finite) simplicial complex $\Delta$ is {\em pure} if each maximal face
of $\Delta$
has the same dimension.
A pure simplicial complex $\Delta$ is {\em Eulerian} if for every
face $F$ of $\Delta$ (including the empty face) the Euler characteristic
of its link
is equal to the Euler characteristic of the sphere of the same dimension,
that is,
$$
\chi(\lk F)=1+(-1)^{\dim(\lk F)}.
$$
In particular,  by 
Poincar\'{e} duality, every odd-dimensional homology manifold is
Eulerian. (Recall that a simplicial complex $\Delta$ is a 
{\em homology manifold} if its
geometric realization $X$ possesses the following property:
for every  $p\in X$ and every $i<\dim X$, $H_i(X,X-p)=0$, while
$H_{\dim X}(X,X-p)\cong \Z$. Here $H_i(X,X-p)$ denotes the $i$-th relative
singular homology with coefficients $\Z$.)

The Upper Bound Conjecture (abbreviated UBC) proposed 
by Motzkin in 1957 (see \cite{Motzkin}) asserts
that among all $d$-dimensional (simplicial) polytopes with $n$ vertices,
the number of $i$-dimensional faces 
(for every $i=1,...,d-1$) is maximized by 
the cyclic polytope $C_d(n)$.
Over the last 40 years this conjecture has been treated extensively by 
many mathematicians:
in 1970, McMullen \cite{MullenUBC} proved
 the UBC for polytopes; McMullen's result was preceded
in 1964 by a surprising work
of Klee, where he verified that the UBC holds for all 
Eulerian complexes with a sufficiently large number of vertices, and
conjectured that it holds for
all Eulerian complexes \cite{Klee1}; in 1975 Stanley proved
the UBC for arbitrary triangulations of spheres
\cite{Stanley1},\cite{Stanley3}, and in 1998 Novik verified
the UBC for triangulations of odd-dimensional
manifolds and several classes of even-dimensional manifolds \cite{Novik}.

In this note we will prove the UBC for a class of odd-dimensional 
simplicial complexes that in particular
includes all odd-dimensional Eulerian  complexes whose geometric
realization has isolated singularities.
More precisely, we obtain the following theorem in which
$f_i(\Delta)$ denotes
the number of $i$-dimensional faces of a complex~$\Delta$, the values
$\beta_i(\Delta)=\dim({\widetilde{H}}_i(\Delta))$ 
denote the
reduced Betti numbers of~$\Delta$ over a field of characteristic 0,
and $C_d(n)$ is a $d$-dimensional cyclic polytope on $n$ vertices.  

\begin{theorem}  \label{UBC-vanish}
Let $\Delta$ be a  pure $(2k+1)$-dimensional simplicial complex on $n$
vertices, such that for
every vertex $v$ of $\Delta$, the link of $v$ is either a
homology manifold whose Euler characteristic is 2,
 or an oriented homology manifold satisfying the following condition
$$
\beta_{k}(\lk v)\leq 2\beta_{k-1}(\lk v)+2\sum_{i=0}^{k-3}\beta_i(\lk v).
$$
Then
$f_i(\Delta) \leq f_i(C_{2k+2}(n))$  for $i=1,2,\ldots, 2k+1$.
\end{theorem}

The main ingredient in the proofs 
is a new invariant of simplicial complexes,
 ${\widetilde{h}}(\Delta)=
({\widetilde{h}}_0, {\widetilde{h}}_1, \ldots, 
{\widetilde{h}}_{\dim(\Delta)})$,
which is a simplicial  analog of the short cubical $h$-vector introduced by
Adin (see \cite{Adin}). 
 We give its definition and list some of its properties in the next
section.
Section \ref{UBC} is devoted to a proof of Theorem~\ref{UBC-vanish}. 
Section~\ref{remarks} contains several remarks and additional results
on the UBC and the $\widetilde{h}$-vector.

\section{$\widetilde{h}$-vector}
     \label{tilde_h_vector}
In this section we introduce the notion of $\widetilde{h}$-vector
for pure simplicial complexes and list some of its properties.
Let us begin by recalling definitions of $f$-vectors and $h$-vectors.
For a $(d-1)$-dimensional simplicial complex $\Delta$, its $f$-vector,
denoted $f(\Delta)$, is a vector 
$(f_{-1}, f_0, f_1, \ldots, f_{d-1})$ where
$f_i$ counts the number of $i$-dimensional faces. 
In particular, 
 $f_{-1}=1$, $f_0$ is the number of vertices of $\Delta$, 
and $f_1$ is the number of edges.
 The $h$-vector of $\Delta$, denoted $h(\Delta)$, is a vector
$(h_0, h_1, \ldots, h_d)$ where
\begin{equation}  \label{h-vector}
h_i(\Delta)=\sum_{j=0}^i (-1)^{i-j}{d-j \choose d-i}f_{j-1}(\Delta), 
    \qquad i=0, 1, \ldots, d.
\end{equation}
Equivalently,
\begin{equation}  \label{f-vector}
f_{j-1}(\Delta)=\sum_{i=0}^j {d-i \choose d-j}h_i(\Delta), 
    \qquad j=0, 1, \ldots, d.
\end{equation}

Ron Adin \cite[eq.~(1), (11)]{Adin}
defined for any cubical complex $C$ its {\em short cubical} h-vector, denoted
$h^{(sc)}(C)=(h^{(sc)}_0, h^{(sc)}_1, \ldots, h^{(sc)}_{\dim(C)}).$
It was later observed by G.~Hetyei that if $C$ is pure, then 
$h^{(sc)}(C)=\sum_{v\in V} h(\lk v)$, 
where $V$ is the set of vertices of $C$.
(Note that the links of the vertices in a cubical complex are simplicial
complexes, and hence the $h$-vector $h(\lk v)$ is well-defined.)

Similarly to the short cubical $h$-vector, we define a short simplicial 
$h$-vector, denoted $\widetilde{h}$, as follows.
\begin{definition}
Let $\Delta$ be a pure $(d-1)$-dimensional simplicial complex
on the vertex set~$V$. 
Define
$$
{\widetilde{h}}(\Delta)=({\widetilde{h}}_0,
{\widetilde{h}}_1, \ldots, {\widetilde{h}}_{d-1}):=\sum_{v\in V} h(\lk v),
$$
so in particular ${\widetilde{h}}_i(\Delta) :=\sum_{v\in V} h_i(\lk v)$.
\end{definition}

The next lemma gives several properties of $\widetilde{h}$.
\begin{lemma}  \label{prop_of_htilde} 
\begin{itemize}
\item[(i)]
Let $\Delta$ be a pure $(d-1)$-dimensional simplicial complex. Then
\begin{eqnarray*}
{\widetilde{h}}_i(\Delta) & = &\sum_{j=0}^i 
   (-1)^{i-j}(j+1){{d-1-j} \choose {d-1-i}}f_j(\Delta) 
\qquad (0\leq i\leq d-1) \;\; \mbox{ and }   \\
 f_j(\Delta) & = & (j+1)^{-1}\sum_{i=0}^j 
   {{d-1-i} \choose {d-1-j}}{\widetilde{h}}_i(\Delta) 
 \qquad (0\leq j\leq d-1).
\end{eqnarray*}
In particular, the $f$-numbers of a simplicial complex
are non-negative linear combinations of its $\widetilde{h}$-numbers.
\item[(ii)]
If $\Delta$ is a pure 
$(2k+1)$-dimensional simplicial complex such that the link of
every vertex is a homology manifold, then the $f$-numbers of $\Delta$
are non-negative linear combinations of 
$\widetilde{h}_0, \widetilde{h}_1, \ldots, \widetilde{h}_{k+1}$.
In other words, 
$$
f_j(\Delta) = \sum_{i=0}^{k+1} b_i^j  \widetilde{h}_i(\Delta) 
                   \qquad 0\leq j \leq 2k+1,
$$
where the coefficients $b_i^j$ are 
independent  of $\Delta$ and are non-negative.
\end{itemize}
\end{lemma}
\proof
Since every $j$-dimensional simplex has $j+1$ vertices,
it follows that
\begin{equation}  \label{f_0i}
\sum_{v\in V} f_{j-1}(\lk v)=(j+1)f_j(\Delta),
\end{equation}
where $V$ is the set of verticies of $\Delta$.
This equation together with
 relations (\ref{h-vector}) and (\ref{f-vector}), 
(applied to the links of vertices) implies part (i).

Part (ii) is a consequence of equation (\ref{f_0i}) and 
\cite[Lemma 6.1]{Novik}, which asserts that the $f$-numbers of
a $2k$-dimensional homology manifold are non-negative linear combinations 
of its $h$-numbers $h_0, h_1, \ldots, h_{k+1}$.
\endproof

\section{The proof of the Upper Bound Theorem}
  \label{UBC}

In this section we prove Theorem \ref{UBC-vanish}.
This will require the following facts and definitions.
\begin{definition}
A simplicial complex $\Delta$ is $l$-neighborly if each set of $l$
of its vertices forms a face in $\Delta$.
\end{definition}

It is  well-known that all $d$-dimensional cyclic polytopes are 
$\lfloor d/2 \rfloor$-neighborly, and that all 
$\lfloor d/2 \rfloor$-neighborly
$d$-dimensional polytopes with $r$ vertices have the same $h$-vector:
$$
h_i=h_{d-i}={r-d+i-1 \choose i} \qquad\mbox{for } 
    0\leq i \leq \lfloor d/2 \rfloor.
$$

In the proof of Theorem \ref{UBC-vanish} we will also use
the following version of the Upper Bound Theorem for even-dimensional
homology manifolds.
\begin{lemma}  \label{h<h}
Let $K$ be a $2k$-dimensional homology manifold on $r$ vertices.
Furthermore, let us assume that either $\chi(K)=2$, or $K$ is
an oriented homology manifold such that
\begin{equation}  \label{beta}
\beta_{k}(K)\leq 2\beta_{k-1}(K)+2\sum_{i=0}^{k-3}\beta_i(K).
\end{equation}
Then 
$$
h_i(K)\leq h_i(C_{2k+1}(r))  \qquad \mbox{ for } 0\leq i \leq k+1.
$$
\end{lemma}
\proof
In the case of $\chi(K)=2$, the lemma follows from \cite[Theorem 6.6]{Novik}
and the Dehn-Sommerville relations for Eulerian complexes \cite{Klee}.
In the second case, the result is a part of the proof 
of \cite[Theorem 6.7]{Novik}.
\endproof

We are now ready to verify Theorem \ref{UBC-vanish}.
The argument is very similar to the proof of a special case
of the cubical upper bound conjecture (see \cite[Theorem 4.3]{BBC}).
The only difference is that we use the $\widetilde{h}$-vector
instead of the short cubical $h$-vector. 

\smallskip\noindent {\it Proof of Theorem \ref{UBC-vanish}: \ }
Let $\Delta$ be a simplicial complex satisfying the conditions of the 
theorem.  By Lemma \ref{prop_of_htilde}(ii), it suffices to check that 
${\widetilde{h}}_i(\Delta) \leq {\widetilde{h}}_i(C_{2k+2}(n))$
for $0\leq i \leq k+1$.  To this end, note that 
for every vertex $v$ of $\Delta$,
$\lk v$ is a simplicial complex on at most $n-1$ vertices
that is either a  homology
manifold with Euler characteristic 2, or
an oriented homology manifold satisfying condition (\ref{beta}).
Thus, by Lemma \ref{h<h}, 
$$
h_i(\lk v)\leq h_i(C_{2k+1}(n-1)) \mbox{ for } 0\leq i \leq k+1.
$$
Since $C_{2k+2}(n)$ is a $(k+1)$-neighborly polytope, it follows that
the link of every vertex of $C_{2k+2}(n)$ is a $k$-neighborly
$(2k+1)$-dimensional polytope on $n-1$ vertices. Hence,
$$
{\widetilde{h}}_i(\Delta)= \sum_{v} h_i(\lk v)\leq 
\sum_{v} h_i(C_{2k+1}(n-1))={\widetilde{h}}_i(C_{2k+2}(n)) \qquad 
 \mbox{ for } 0\leq i \leq k+1,
$$
implying the theorem.
\endproof

\begin{corollary}
Let $\Delta$ be a $(2k+1)$-dimensional oriented pseudomanifold on
$n$ vertices
such that the link of every vertex is either a $2k$-dimensional homology
manifold  with vanishing middle homology, or it is
a $2k$-dimensional homology manifold whose Euler characteristic 
$\chi$ satisfies $(-1)^k (\chi-2)\leq 0$.
Then
$$
f_i(\Delta)\leq f_i(C_{2k+2}(n)) \qquad \mbox{for } 1\leq i\leq 2k+1.
$$
\end{corollary}
\proof
Any such complex $\Delta$ satisfies the 
assumptions of Theorem \ref{UBC-vanish}.
\endproof

\section{Additional remarks and results}   \label{remarks}
\paragraph{\bf 1.}
Theorem \ref{UBC-vanish} proves a special case of Gil
Kalai's conjecture \cite[Section 7]{Novik}
that the UBC holds for all simplicial
complexes having the property that every link (of a face)
of dimension $2k$ ($k=1,2,\ldots$)
satisfies condition (\ref{beta}).

\medskip
\paragraph{\bf 2.}
In his proof of the UBC for spheres \cite{Stanley1}, \cite{Stanley3}, 
R.~Stanley showed that if $K$
is a $(d-1)$-dimensional homology sphere on $n$ vertices then
\begin{equation}   \label{sphere}
h_i(K)\leq h_i(C_d(n)) \qquad\mbox{for } 0\leq i\leq d-1.
\end{equation}
Since the $f$-numbers of any simplicial complex $\Delta$ are non-negative
combinations of its $\widetilde{h}$-numbers (by Lemma \ref{prop_of_htilde}(i)),
arguing exactly  as in the proof of Theorem~\ref{UBC-vanish}, but using
(\ref{sphere}) instead of Lemma \ref{h<h}, we obtain a new proof of the UBC
for odd-dimensional homology manifolds. 
This proof is shorter and more elementary
than the one presented in \cite[Theorem 1.4]{Novik}. (It does not
use any facts about Buchsbaum complexes!)

\medskip
\paragraph{\bf 3.} 
It would be interesting to clarify whether for a $(2k+1)$-dimensional
complex~$\Delta$
satisfying the assumptions of Theorem \ref{UBC-vanish},
the inequality $h_i(\Delta)\leq h_i(C_{2k+2}(n))$ ($0\leq i\leq k+1$)
necessarily holds.  We have the expression
\begin{eqnarray*}
h_r(\Delta) &=&\sum_{j=0}^r(-1)^{r-j}{2k+2-j \choose 2k+2-r}f_{j-1}(\Delta) \\
&=& (-1)^r{2k+2 \choose r}+
        \sum_{i=0}^{r-1}\widetilde{h_i}(\Delta){2k+1-i\choose 2k+2-r}
	\sum_{j=i+1}^r \frac{1}{j}(-1)^{r-j}{r-i-1\choose r-j}\\
&=& (-1)^r{2k+2 \choose r}+
   \sum_{i=0}^{r-1}\widetilde{h_i}(\Delta){2k+1-i\choose 2k+2-r}
	  \int_0^1 x^{i} (x-1)^{r-i-1} \,dx  .\\
\end{eqnarray*}
Hence the coefficients of $\widetilde{h}$-numbers in the
expression for $h_r$ alternate in sign
so that short simplicial $h$-vectors are not sufficient to resolve this
question.

\medskip
\paragraph{\bf 4.} {\bf Lower bounds.}
Let $\Delta$ be a simplicial complex, let $\skel_i(\Delta)$ denote
its $i$-dimensional skeleton,
and let $\chi_i(\Delta):=\chi(\skel_i(\Delta))=\sum_{j=0}^i (-1)^j f_j(\Delta)$
denote the Euler characteristic of $\skel_i(\Delta)$.
It was shown in \cite{Novik2} that if $\Delta$ is a $(2k-1)$-dimensonal
manifold, then $(-1)^i\chi_i(\Delta)\geq 0$ for $0\leq i\leq 2k-1$.
The proof relied on several facts about Buchsbaum complexes.
Using $\widetilde{h}$-numbers we provide a short proof of the following
related result.
\begin{proposition}
Let $\Delta$ be a $(d-1)$-dimensional Buchsbaum simplicial complex
(i.e. a pure simplicial complex such that for every vertex $v\in\Delta$
the link of $v$ is Cohen-Macaulay). Then
$(-1)^i\chi_i(\Delta)\geq 0$ for $0\leq i\leq \lfloor(d-1)/2\rfloor$.
\end{proposition}
\proof
Since for every vertex $v\in\Delta$, $\lk v$ is Cohen-Macaulay, 
it follows 
that $h_i(\lk v)\geq 0$ for $i=0,1,\ldots,d-1$, and hence,
${\widetilde{h}}_i(\Delta)\geq 0$ for $i=0,1,\ldots,d-1$.
Expressing the $f$-numbers of $\Delta$ in terms of its
$\widetilde{h}$-numbers (Lemma \ref{prop_of_htilde}(i)), we obtain
\begin{equation}  \label{pos}
(-1)^i\chi_i(\Delta)=\sum_{j=0}^i(-1)^{i-j}f_j=
\sum_{l=0}^i \left( \sum_{j=l}^i (-1)^{i-j}\frac{1}{j+1} 
{d-1-l \choose d-1-j}\right){\widetilde{h}}_l.
\end{equation}
It is straightforward to show that if $0\leq i\leq \lfloor(d-1)/2\rfloor$
and $0\leq l\leq i$,
then
$$
\frac{1}{i+1}{d-1-l \choose d-1-i} \geq \frac{1}{i}{d-1-l \choose d-i}\geq
\ldots\geq \frac{1}{l+1}{d-1-l \choose d-1-l}.
$$
Hence for any $0\leq i\leq \lfloor(d-1)/2\rfloor$,
all coefficients of $\widetilde{h}$-numbers in equation (\ref{pos})
are non-negative, implying the proposition.
\endproof

\medskip
\paragraph{\bf 5.} {\bf Semi-Eulerian complexes.}
One may also use short simplicial 
$h$-vectors and the Dehn-Sommerville relations
to give a new proof of the fact that all
odd-dimensional semi-Eulerian simplicial (or regular cell)
complexes are Eulerian.  This result was proven more generally
for posets in \cite[Exercise 3.69(c)]{Stanley} by a very
different approach.

\section*{Acknowledgments}
We would like to thank Victor Klee and Richard Stanley for their
comments on the preliminary version of this paper. We are also
grateful to Richard Stanley for directing us to the topic of cubical
complexes.

\end{document}